\documentclass[amssymb,12pt ]{article}
\usepackage{amscd}
\usepackage{amsfonts}
\usepackage{amssymb}
\usepackage{amsmath}
\usepackage{geometry} 
\usepackage{graphicx}

\usepackage{longtable}

\newtheorem{thm}{Theorem}[section]

\newtheorem{rem}[thm]{Remark}

\newenvironment{proof}{{\bf Proof }}{\hfill\framebox[2mm]{}}

\begin{document}

\title{Volume of $C^{1,\alpha}$-boundary domain in extended hyperbolic space}
\author{Yunhi Cho\footnote{The first author was supported by Research Fund 2004 of University of Seoul.} and
Hyuk Kim\footnote{The second author was supported by grant no.(R01-1999-000-00002-0)
  from the Basic Research Program of the Korea Science $\&$ Engineering Foundation.}}
\date{}
\maketitle

\begin{abstract}
We consider the projectivization of Minkowski space with the
analytic continuation of the hyperbolic metric and call this an
extended hyperbolic space. We can measure the volume of a domain
lying across the boundary of the hyperbolic space using an
analytic continuation argument. In this paper we show this method
can be further generalized to find the volume of a domain with
smooth boundary with suitable regularity in dimension 2 and 3. We
also discuss that this volume is invariant under the group of
hyperbolic isometries and that this regularity condition is sharp.
\end{abstract}

\renewcommand \figurename{Fig.}
\makeatletter
\renewcommand\fnum@figure[1]{\textbf{\figurename} }
\makeatother

\section{Introduction and preliminaries}

In [1] we considered an extended model of hyperbolic space and
studied how we can define a volume of a domain which lies beyond
the infinity of the hyperbolic space. Such investigation gives us
a natural way of studying various geometric objects in Lorentz
geometry in a manner consistent with those in hyperbolic geometry.
The method of calculating volume of such domain is essentially an
analytic continuation argument and works very well with a domain
with analytic boundary. But if the boundary is smooth or just
continuous, then the volume problem turns out to be very delicate
and the required regularity of the boundary necessary for
finiteness of volume depends on the dimension. We discuss this
phenomenon in detail in this paper focusing especially on
dimension two or three. Then we discuss the invariance of the
volume of domains which has boundary with necessary regularity in
these dimensions. We keep the same notations used in [1], but we
provide necessary materials in detail so that the paper is as
self-contained as possible and can be read independently from [1].
And here we do not intend to mention why the extended model is
natural and what applications we can obtain using this model. We
refer the reader to the paper [1] for all these explanations and
other references as well.

Let $\Bbb R^{n,1}$ denote the Minkowski space, i.e., $\Bbb
R^{n+1}$ with the inner product of signature $(n,1)$ given by
\[\langle x,y\rangle=-x_0y_0+x_1y_1+\cdots +x_n y_n.\]
The hyperbolic space, Lorentz space and the light cone are defined
as the sets $\{x\in \Bbb R^{n,1}|\langle x,x\rangle=\alpha\}$ with
$\alpha=-1,1,0$ respectively together with the induced metric. If
we project these sets radially to an affine subspace $\Bbb
K^n:=\{1\}\times\Bbb R^n \subset\Bbb R^{n,1}$, then we obtain a
unit ball as Kleinian projective model for hyperbolic space $\Bbb
H^n$, Lorentz space of constant sectional curvature 1 outside the
ball and the light cone as the common boundary $\partial \Bbb H^n$
of these two spaces.

\begin{figure}[h]
\begin{center}
\includegraphics[width=0.35\textwidth]{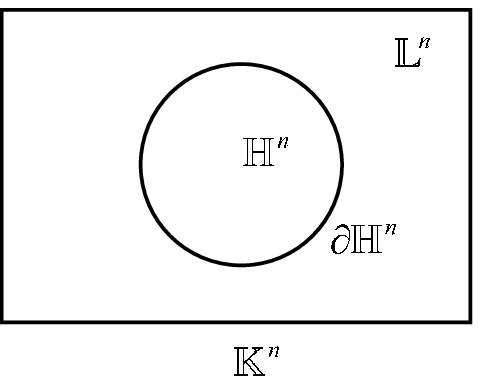}
\caption{\textbf{1}}
\end{center}
\end{figure}

If we change the sign of the induced metric on the Lorentz space,
then the new Lorentz space denoted by $\Bbb L^n$ has constant
sectional curvature $-1$ and the metrics on both parts $\Bbb H^n$
and $\Bbb L^n$ have the exactly same formula on $\Bbb K^n$
\[ds^2_K=\left({\Sigma x_i dx_i\over 1-|x|^2}\right)^2 +{\Sigma dx_i^2\over 1-|x|^2}.\]
And the induced volume form is given by
\[dV_K=\frac{dx_1\wedge\cdots\wedge
dx_n}{(1-|x|^2)^{\frac{n+1}{2}}}.\]

Now for a domain $U$ in $\Bbb H^n$, the volume of $U$ will be
simply given by the integration of $dV_K$ on $U$. For a domain $U$
lying across the boundary of $\Bbb H^n$, we formally calculate the
volume of $U$ using the polar coordinates as follows:
$$
\aligned \text{vol} (U)&=\int_U  \frac{dx_1\cdots dx_n}{(1-|x|^2)^{\frac{n+1}{2}}} \\
              &=\int_{G^{-1}(U)} \frac{r^{n-1}}{(1-r^2)^{\frac{n+1}{2}}} dr
              d\theta\\
&=\int_a^b \frac{r^{n-1}F(r)}{(1-r^2)^{\frac{n+1}{2}}} dr,\quad
F(r)=\int_{G^{-1}(U)\cap S^{n-1}(r)} d\theta,\\
\endaligned
$$
where $G:(r,\theta)\mapsto(x_1,\cdots ,x_n)$ is the  polar
coordinates, $S^{n-1}(r)$ is the Euclidean sphere of radius $r$
and $d\theta$ is the volume form of the Euclidean unit sphere $
\Bbb S^{n-1}$.

 Now this integral with
respect to $r$ does not make sense in general, but for a domain
$U$ with analytic boundary transversal to $\partial\Bbb H^n$ we
may use contour integral to define a volume of $U$.
$$
\text{vol} (U):= \int_{\gamma}
\frac{r^{n-1}F(r)}{(1-r^2)^{\frac{n+1}{2}}} dr,
$$
where $\gamma$ is a contour given by
$$
 \gamma(t)=\left\{\aligned
                   &t,\qquad \qquad \quad a\le t\le 1-\delta,\\
                   &1+\delta e^{\frac{i(1-t)\pi}{\delta}},  1-\delta\le t\le 1,\\
                   &t+\delta, \qquad\quad\ 1\le t\le b-\delta.\\
                   \endaligned \right.
$$
\begin{figure}[h]
\begin{center}
\includegraphics[width=0.45\textwidth]{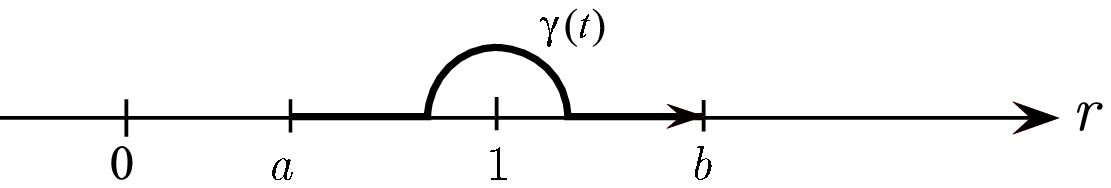}
\caption{\textbf{2}}
\end{center}
\end{figure}

 Note that the analyticity and transversality of the
boundary of $U$ was needed to make sure $F(r)$ is an analytic
function of $r$ near $r=1$. For a domain $U$ in the Lorentz part,
our choice of the contour $\gamma$ naturally determines the sign
of vol($U$) as $i^{n+1}$ and so is determined the sign of $dV_K$
(see [1]).

In [1], it is shown that vol($U$) can also be obtained through a
complex approximation. Let
$$
ds^2_{\epsilon}=\left({\Sigma x_i dx_i\over
d^2_{\epsilon}-|x|^2}\right)^2 +{\Sigma dx_i^2\over
d^2_{\epsilon}-|x|^2},
$$
where $d_{\epsilon}=1-{\epsilon}i$ with $\epsilon>0$ and
$i=\sqrt{-1}$, so that $ds_K^2=\lim_{\epsilon\to
0}ds^2_{\epsilon}$. Then the induced volume form is given by
$$
dV_{\epsilon}=\frac{d_{\epsilon}dx_1\wedge\cdots\wedge dx_n}
{(d^2_{\epsilon}-|x|^2)^{\frac{n+1}{2}}}
$$
and let $\mu(U):=\lim_{\epsilon\to 0}\int_U dV_{\epsilon}$. Here
the choice of sign of $dV_{\epsilon}$ is determined by the
continuity on $\epsilon\ge 0$ and the sign of $dV_K$. Then it was
shown in [1, Proposition 2.1 and 3.2] that $\mu$ is finitely
additive and $\mu(U)=\text{vol }(U)$ for a domain $U$ with an
analytic boundary transversal to $\partial\Bbb H^n$. We actually
show this fact in the next section in a different model. The
finite additivity follows easily from the definition of $\mu$.
Also notice that if $U$ is a domain lying solely in $\Bbb H^n$ or
$\Bbb L^n$, then
$$
\mu(U)=\lim_{\epsilon\to 0}\int_U
dV_{\epsilon}=\int_U\lim_{\epsilon\to 0} dV_{\epsilon}=\int_U dV_K
$$
by the Lebesgue dominated convergence theorem and coincides with
the usual volume.

 The measure theory for $\mu$ seems to be very
delicate and it is not easy to find a large enough class of
$\mu$-measurable sets, that is, Lebesgue measurable sets with
$\mu(U)<\infty$. The present work reflects the effort of finding
and explaining more about $\mu$-measurable sets and we find that
a domain with $C^{1,\alpha}$ boundary in dimension 3 ($C^{0,\frac
12 +\alpha}$ boundary for dimension 2, respectively) is actually
$\mu$-measurable, and also show that this regularity condition is
in fact sharp.

\section{A flattened model for computation}
We prove the results stated in the previous section by computing
various integrals. But computing the integral whose singularities
lies on the unit sphere in $\Bbb K^n$ is certainly inconvenient
and we want to introduce a new model to facilitate the
computation. In this model, we want the singularity sets of our
volume form is a hyperplane. The immediate choice is a Cayley
transformation or a reflection $\sigma$ with respect to a sphere
of radius $\sqrt{2}$ with the center at
$e_n=(0,\ldots,0,1)\in\Bbb K^n$.

We see immediately that under the reflection $\sigma$, $\Bbb H^n$
is sent to the lower half space and $\Bbb L^n$ to the upper half
space.

From the obvious identities,
$$
\left\{\aligned
                   &y-e_n=\lambda(x-e_n),\quad\lambda\in\Bbb R,\\
                   &|y-e_n||x-e_n|=2,\\
                   \endaligned \right.
$$

\begin{figure}[h]
\begin{center}
\includegraphics[width=0.6\textwidth]{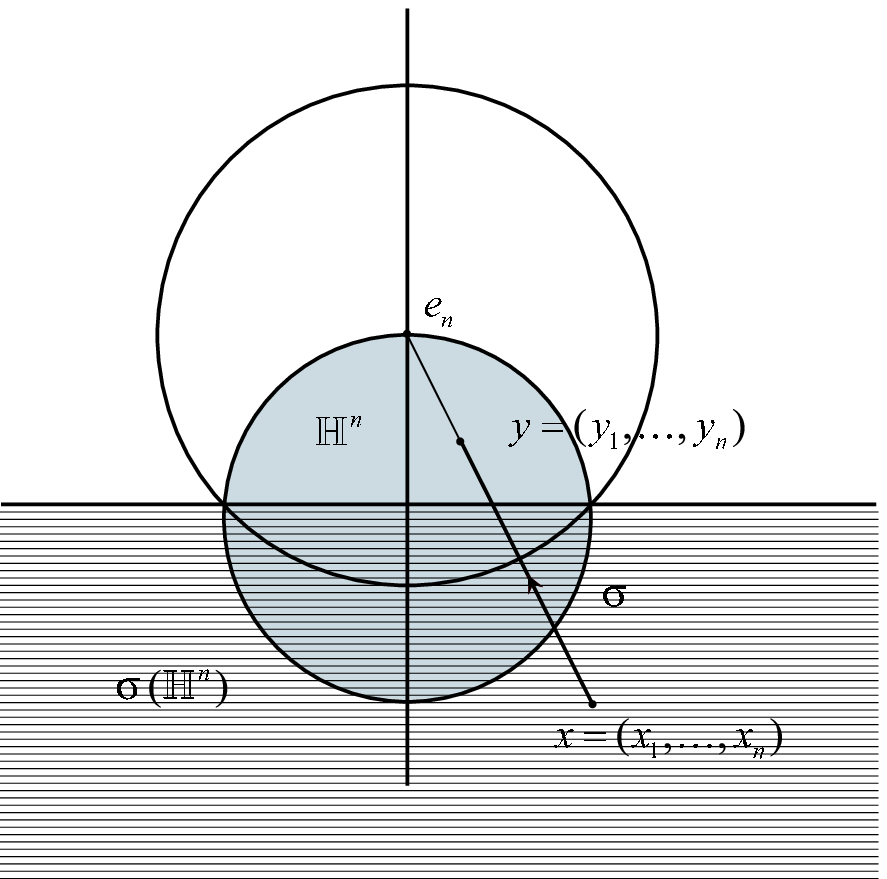}
\caption{\textbf{3}}
\end{center}
\end{figure}

\noindent{we easily obtain that $y=\sigma(x)$ is given by}
$$
\sigma : \left\{\aligned
                   &y_i = {2x_i \over |x-e_n|^2},\quad i=1,\ldots,n-1,\\
                   &y_n ={2(x_n-1) \over |x-e_n|^2} +1.\\
                   \endaligned \right.
$$
We compute directly using this formula that the metric $ds_K^2$
is pulled back by $\sigma$ to
\[ds^2=\sigma^*(ds^2_K)=\left({\alpha dx_n- x_n d\alpha\over 2\alpha x_n}\right)^2 -{\Sigma dx_i^2\over \alpha x_n},\]
where $\alpha=|x-e_n|^2=x_1^2+\cdots+x_{n-1}^2+(x_n-1)^2$ so that $d\alpha=2(\Sigma x_i dx_i -dx_n).$

Also the volume form $dV_K$ is pulled back to
\[dV=\sigma^*(dV_K)=-\frac{dx_1\wedge\cdots\wedge dx_n}{2(-x_n)^{\frac{n+1}{2}}\alpha^{\frac{n-1}{2}}}.\]
Here notice that the first negative sign appears since $\sigma$
is orientation reversing and we can ignore  this when we compute
the integrals for volume. If $x_n>0$, that is, if $x\in\Bbb L^n$,
we need to determine the sign of $(-1)^{\frac{n+1}{2}}$, and this
should be determined as $(-i)^{n+1}$ in order to give the sign of $dV$ as $i^{n+1}$ as given in the previous section.

This new model $\Bbb E^n$ is of course quite different from the
Poincar\'e half space model. It is clear from the construction
that the geodesics in this model are the circles (including lines
viewed as a special case of circles passing through the infinity)
passing through the point $e_n$, and more generally spheres
(including planes) passing through $e_n$ are the totally geodesic
submanifolds.

Let's consider first the volume of a domain $U$ with analytic
boundary transversal to $\partial\Bbb H^n$ in the new model $\Bbb
E^n$. Note that $\sigma^{-1}=\sigma$.
$$\mu(U)=\lim_{\epsilon\to 0}\int_{\sigma(U)}dV_{\epsilon}=\lim_{\epsilon\to 0}\int_U d\tilde V_{\epsilon},$$
where
$$d\tilde V_{\epsilon}=\sigma^*(dV_{\epsilon})=-\frac{1-\epsilon i}{2}\frac{dx_1\wedge\cdots\wedge dx_n}
{(\frac{-\epsilon^2-2\epsilon
i}{4}\alpha-x_n)^{\frac{n+1}{2}}\alpha^{\frac{n-1}{2}}}.$$ We also ignore negative sign in the above
formula of $d\tilde V_{\epsilon}$ when we compute integrals. The
induced volume form  $d\tilde V_{\epsilon}$ has a complicated
formula, and instead we use a different  simple volume
approximation $d\mu_{\epsilon}$ which gives us the same
$\mu$-measure of $U$.

\begin{thm}\label{a} Let $U$ be a bounded domain with analytic boundary transversal
to $\partial\Bbb H^n$ in $\Bbb E^n$ and let
$$
d\mu_{\epsilon}=\frac{dx_1\wedge\cdots\wedge dx_n}{2(-x_n-\epsilon i)^{\frac{n+1}{2}}\alpha^{\frac{n-1}{2}}},
\quad \alpha=|x-e_n|^2.
$$
Then $$\mu(U)=\lim_{\epsilon\to 0}\int_U d\mu_{\epsilon}.$$
Furthermore for the domain $U$ with $-\delta<x_n<\delta$,
$$
\mu(U)=\int_{\gamma}\int \frac{dx_1\wedge\cdots\wedge
dx_{n-1}}{2(-x_n)^{\frac{n+1}{2}}\alpha^{\frac{n-1}{2}}}~ dx_n,
$$
where $\gamma$ is a contour given below in Fig. 4.
\end{thm}

\begin{proof}
First observe that the volume of a domain lying completely inside
of $\Bbb H^n$ or $\Bbb L^n$, the same statement holds. This can be
easily checked from Lebesgue dominated convergence theorem using
$|x_n+\epsilon i|\ge |x_n|$ and from that $d\mu_0$ is just $dV$.
Now by the finite additivity of the volume, we may assume that $U$
lies in the domain $-\delta<x_n<\delta$ for a sufficiently small
$\delta>0$. We will prove the theorem in the following two steps:

{\bf Step 1:} {\it $$\mu(U):=\lim_{\epsilon\to 0}\int_U\tilde
V_{\epsilon}=\lim_{\epsilon\to 0}\int_{-\delta}^{\delta}\!\int
 d\tilde V_{\epsilon}=\lim_{\epsilon\to 0}\int_{\gamma}\int d\tilde V_{\epsilon}$$
and
$$\lim_{\epsilon\to 0}\int_U d\mu_{\epsilon}=\lim_{\epsilon\to 0}\int_{-\delta}^{\delta}\!\int
d\mu_{\epsilon}=\lim_{\epsilon\to 0}\int_{\gamma}\int
d\mu_{\epsilon}.$$ Here for the double integral
$\int_{-\delta}^{\delta}\!\int $ we integrate first with respect
to the variables $(x_1,\ldots,x_{n-1})$ and then with respect to
the variable $x_n$.}

The contour integral  $\int_{\gamma}$ is an integration with
respect to complex variable $x_n$ and $\gamma$ is a contour given
below in Fig. 4.

\begin{figure}[h]
\begin{center}
\includegraphics[width=0.45\textwidth]{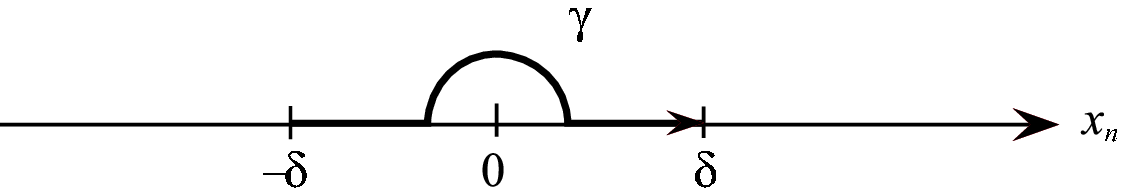}
\caption{\textbf{4}}
\end{center}
\end{figure}

 {\bf Step 2:} {\it $$\lim_{\epsilon\to
0}\int_{\gamma}\int d\tilde V_{\epsilon}=\int_{\gamma}\int
d\tilde V_0$$ and
$$\lim_{\epsilon\to 0}\int_{\gamma}\int d\mu_{\epsilon}=\int_{\gamma}\int
d\mu_0.$$}

Note that
$$d\tilde V_0:=\lim_{\epsilon\to 0}d\tilde V_{\epsilon}=dV=\lim_{\epsilon\to 0} d\mu_{\epsilon}=:d\mu_0,$$
and hence the theorem follows from Step 1 and Step 2.

{\bf Proof of Step 1:} We can show that
$\int_{-\delta}^{\delta}\!\int d\tilde
V_{\epsilon}=\int_{\gamma}\int d\tilde V_{\epsilon}$ if we can
show that the pole of $\int d\tilde V_{\epsilon}$ as a function of
$x_n$ has negative imaginary part for all $\epsilon>0$. This looks
intuitively so because
 $$
dV_{\epsilon}=\frac{d_{\epsilon}r^{n-1}dr d\theta_1\cdots d\theta_{n-1}} {(d_{\epsilon}^2-r^2)^{\frac{n+1}{2}}}
$$
in spherical coordinates has pole with negative imaginary part
near $r=1$ and $r$ corresponds essentially to $x_n$ under the
coordinate change map which is real.

To be more precise, let $g$ be the coordinate change map
$x=(x_1,\ldots,x_n)=g(r,\theta_1,\ldots,\theta_{n-1})$ given by a
composite of spherical coordinates and the reflection $\sigma$:
$$
\left\{\aligned
y_1&=r \sin\theta_1\cdots \sin\theta_{n-2}\sin\theta_{n-1},\\
y_2&=r \sin\theta_1\cdots \sin\theta_{n-2}\cos\theta_{n-1},\\
\vdots\\
y_{n-1}&=r \sin\theta_1\cos\theta_2,\\
y_n&=r \cos\theta_1,
\endaligned\right.
$$
with $r>0$, and
$$
x_1=\frac {2}{\alpha}~y_1, \cdots,x_{n-1}=\frac {2}{\alpha}~y_{n-1}, x_n=\frac {2}{\alpha}(y_n-1)+1,
$$
with $\alpha=|y-e_n|^2=y_1^2+\cdots +y_{n-1}^2+(y_n-1)^2.$

Write
$$d\tilde V_{\epsilon}=\frac{1}{f_{\epsilon}(x_1,\ldots,x_n)}dx_1\wedge\cdots\wedge dx_n$$
 and consider zeroes of $f_{\epsilon}(c_1,\ldots,c_{n-1},x_n),~~c_i\in\Bbb R.$
 We claim that  $f_{\epsilon}(c_1,\ldots,$ $c_{n-1},x_n)$ has no real zeroes. Indeed if it had,
 $f_{\epsilon}\circ g$ would have real zeroes since $g$ is real and hence
$dV_{\epsilon}=\frac{d_{\epsilon}r^{n-1}dr d\theta}
{(d_{\epsilon}^2-r^2)^{\frac{n+1}{2}}}=\frac{1}{f_{\epsilon}}\circ
g(\text{det}~g')~dr d\theta$ would have real poles, which is
absurd.

Therefore the imaginary part of zeroes of
$f_{\epsilon}(c_1,\ldots,c_{n-1},x_n)$ is either positive or
negative on a connected open set consisting of parameters
$(c_1,\ldots,c_{n-1})$ by continuity, and we can determine the
sign by checking at one point. Notice that the $r$-axis given by
$\theta_1=\pi$ is sent to $x_n$-axis $(x_1=\cdots=x_{n-1}=0)$
under $g$. In fact, $x_n=\frac{r-1}{r+1},~r>0$ and this is an
increasing function of $r$. If we complexify the real analytic
function $g$, the complex analytic function $g_{\Bbb C}$ will
preserve the negative imaginary parts and send $\{\text{im}
~r<0\}$ to $\{\text{im} ~x_n<0\}$ by the orientation reasoning.
In this argument, the point $(0,\ldots,0,x_n)$ does not belong to
the natural domain, i.e., the image under $g$ of a maximal
connected open domain where $g$ is 1-1, but it is a boundary
point of such domain, and the negativity of imaginary part of
zeroes still follows.

Now since $d\tilde
V_{\epsilon}=\frac{1}{f_{\epsilon}(x_1,\ldots,x_n)}dx_1\wedge\cdots\wedge
dx_n$ has  poles with negative imaginary part for all
$x_1=c_1,\cdots,x_{n-1}=c_{n-1},~c_i\in \Bbb R$, therefore $\int
\frac{1}{f_{\epsilon}(x_1,\ldots,x_n)}dx_1\wedge\cdots\wedge
dx_{n-1}$ as a function of $x_n$ is analytic near $x_n=0$ with the
poles only in the negative imaginary part. Here the analyticity
comes from the condition that $U$ has an analytic boundary
transversal to $\partial\Bbb H^n$.

The proof of (2) is immediate by the same pole argument.

{\bf Proof of Step 2:} For this part, we use Lebesgue dominated
convergence theorem and it suffices to show when $U$ is a compact
set, say $U=D\times\gamma\subset\Bbb R^{n-1}\times\Bbb C$ with
$D$ compact domain. We essentially are integrating on a domain
near $r=1$ in Kleinian model $\Bbb K^n$ which is symmetric with
respect to the rotations and hence we may assume the corresponding
domain $U$ in $\Bbb E^n$ is a small compact set near $x_n=0$
without loss of generality. Since
$dV_{\epsilon}=\frac{d_{\epsilon}r^{n-1}dr d\theta}
{(d_{\epsilon}^2-r^2)^{\frac{n+1}{2}}}$ is clearly uniformly
bounded (with respect to $\epsilon>0$) on $g^{-1}(U)$, its pull
back $d\tilde V_{\epsilon}=g^{-1*}d V_{\epsilon}$, only differing
by Jacobian determinant, is also uniformly bounded. Hence
Lebesgue dominated convergence theorem applies.

The proof of (2) is clear by the same argument.
\end{proof}

\begin{rem}\label{f}
The boundedness condition for a domain $U$ in the statement of
Theorem 2.1 is rather superficial. For a domain in the extended
hyperbolic space, the finiteness of the volume depends only on
how it crosses the boundary of $\Bbb H^n$. And by the finite
additivity of $\mu$, it suffices to consider a small domain near
$\partial\Bbb H^n$ in $\Bbb K^n$, which we may assume is bounded
in $\Bbb E^n$ by considering rotation in $\Bbb K^n$ if necessary
before applying Cayley transformation to $\Bbb E^n$.
\end{rem}

\begin{rem}\label{g}
The analyticity is required only to guarantee that the integral
first with respect to variables $(x_1,\ldots,x_{n-1})$ viewed as a
function of $x_n$ is analytic to replace the second integral by
the contour integral. Therefore as far as this first integral on a
domain $U$ is an analytic function of $x_n$, the proof works
equally well. In fact, Theorem 2.1 can be generalized to the case
when $U$ has a piecewise analytic boundary transversal to
$\partial\Bbb H^n$ (see [1, Proposition 3.2]).
\end{rem}

\section{Volume of a domain with $C^{1,\alpha}$-boundary}

In this section we want to show that a domain $U$ passing through
$\partial\Bbb H^n$ with suitable regularity has a finite volume,
i.e., $\mu(U)<\infty$ by computing in the flattened model. We
first consider the case of dimension 2 and then the more
complicated case of dimension 3.

A domain $U$ in $\Bbb K^n$ will be said to be {\it
$C^{k,\alpha}$-transversal to $\partial\Bbb H^n$} if the boundary
of $U$ is given locally as a $C^{k,\alpha}$ function near
$\partial\Bbb H^n$ and transversal to $\partial\Bbb H^n$ in the
usual sense if $k\ge 1$. Namely for each point $p$ in the
intersection of $U^b$, the boundary of $U$, and $\partial\Bbb
H^n$, there is a neighborhood $V$ of $p$ such that $U^b\cap V$
can be written as a zero set of a single $C^{k,\alpha}$-function
which is transversal to $\partial\Bbb H^n$. In the case of
dimension 2, a domain $U$ in $\Bbb K^2$ is {\it
$C^{0,\alpha}$-transversal to $\partial\Bbb H^2$} if locally the
boundary of $U$ near $\partial\Bbb H^2$ can be written as
$\theta=g(r)$ for a $C^{0,\alpha}$ function $g$.

In the following discussions, we will say for the sake of
convenience that an integral $\int f$ is equivalent to $\int g$,
denoted by $\int f\sim \int g$, if $\int f<\infty$ holds iff $\int
g<\infty$.

\begin{thm}\label{b} In the two-dimensional extended hyperbolic space, the area of a domain $U$
which is $C^{0,\frac 12 +\alpha}$-transversal to $\partial\Bbb
H^2$ is finite.
\end{thm}

\begin{proof}
We will compute in the flattened model and transversality
condition for the boundary may be written as $x_1=g(x_2)$ for a
$C^{0,\frac 12 +\alpha}$ function $g$. It suffices to consider
the $C^{0,\frac 12 +\alpha}$-transversal domain $U$ in the
flattened model which can be divided into pieces, one parallel
strip perpendicular to $x_1$-axis and other pieces (at most four
pieces) with only one vertex lying in $x_1$-axis.

\begin{figure}[h]
\begin{center}
\includegraphics[width=0.55\textwidth]{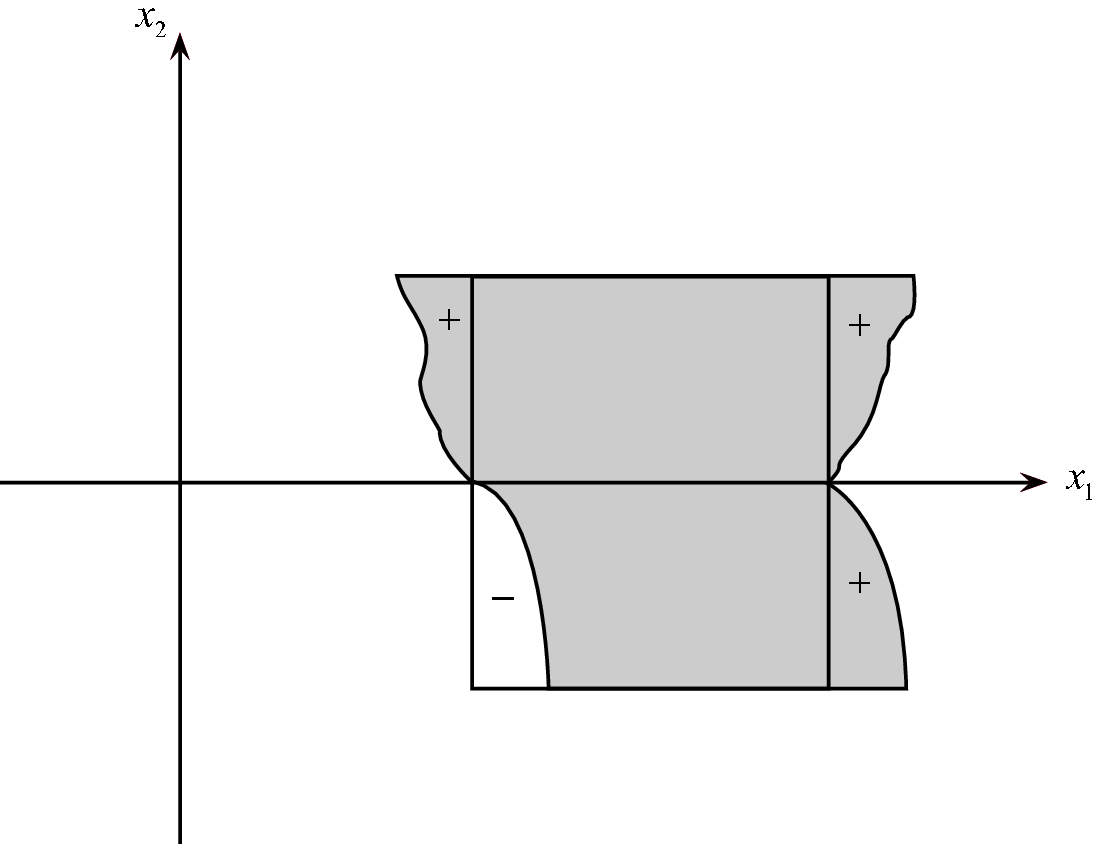}
\caption{\textbf{5}}
\end{center}
\end{figure}

 The transversal strip has finite area by the
Theorem \ref{a} and it suffices to show that each one vertex
domain has also finite area. This can be done if we can show for
a function $x_1=g(x_2)$ with $g(0)=0$ that
$$
\int_0^{\delta}\!\!\int_0^{g(x_2)}\frac{dx_1}{x_2^{\frac
32}(x_1^2+(x_2-1)^2)^{\frac 12}}~dx_2<\infty .
$$
This integral is clearly equivalent to
$\int_0^{\delta}\frac{g(x_2)}{x_2^{\frac 32}}~dx_2$ and by
$C^{0,\frac 12 +\alpha}$ condition of $g(x_2)$, we have
$|g(x_2)|\le C|x_2|^{\frac 12 +\alpha}$ and hence
$$\int_0^{\delta}\frac{|g(x_2)|}{x_2^{\frac 32}}~dx_2\le
C\int_0^{\delta}\frac{1}{x_2^{1-\alpha}}~dx_2<\infty.$$ Thus
$\int_0^{\delta}\frac{g(x_2)}{x_2^{\frac 32}}~dx_2$ is finite.
\end{proof}

So every polygonal domain transversal to $\partial\Bbb H^2$ has
finite area trivially.

\begin{thm}\label{b} In the three-dimensional extended hyperbolic space, the volume of a domain $U$
which is $C^{1,\alpha}$-transversal to $\partial\Bbb H^3$ is
finite.
\end{thm}

\begin{proof}
We work in the flattened model as before. We first explain our
strategy for the proof schematically in dimension 2 since the
three dimensional picture is more complicated. Since all the
difficulties arise near the boundary and near the hyperplane
$\partial\Bbb H^n=\{x|x_n=0\}$, we first localize the problem by
taking a small rectangle near boundary in $\partial\Bbb
H^n=\{x|x_n=0\}$ and we want to show that the volume of the
shaded domain in following picture is finite. We prove this by
showing each of the following three types of integrals ((1), (2)
and (3) in Fig. 6) have finite values.

\begin{figure}[h]
\begin{center}
\includegraphics[width=0.7\textwidth]{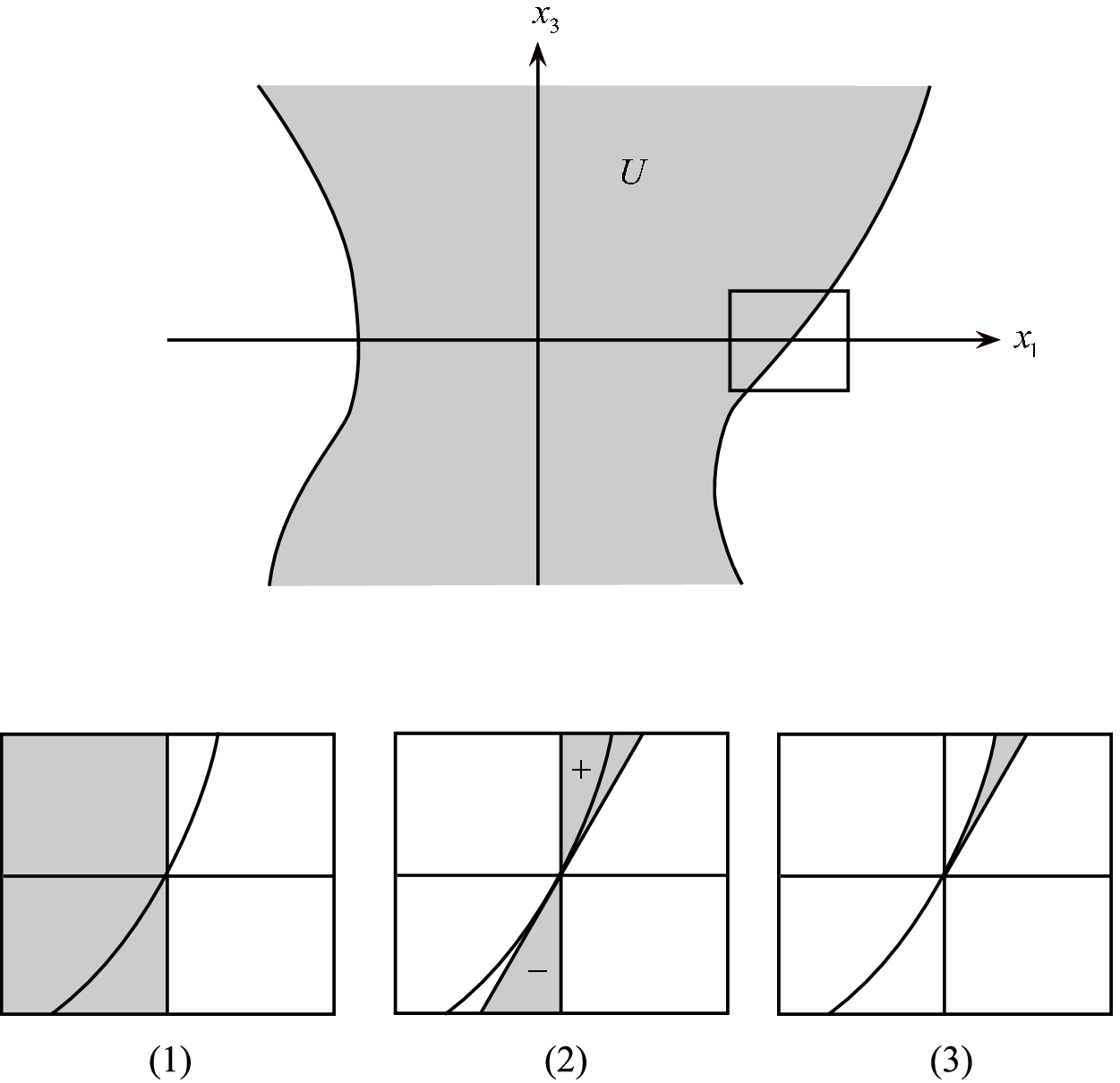}
\caption{\textbf{6}}
\end{center}
\end{figure}

In dimension 3, the basic idea is the same as above and we use
boxes instead of rectangles. We still have the above three types
of integrals and show each of three is finite. We have to prove
the following three integrals (1), (2) and (3) are finite.

For a compact domain $B$ with $C^1$ boundary in the plane
$x_3=0$, the volume of a cylindrical domain $B\times
[-\delta,\delta]$ is represented by $\lim_{\epsilon \to
0}\int^{\delta}_{-\delta}\int_B~d\tilde V_{\epsilon}$, and

\begin{equation}\label{1} \lim_{\epsilon \to 0}\int_B\int^{\delta}_{-\delta}~d\tilde
V_{\epsilon}=\lim_{\epsilon \to 0}\int_B\int_{\gamma}~d\tilde
V_{\epsilon}=\int_B\int_{\gamma}~d\mu_0~<\infty.
\end{equation}
This follows from the pole argument used in the proof of Theorem
2.1 and uniform boundedness of $F_{\epsilon}$ on compact set,
where $d\tilde
V_{\epsilon}=F_{\epsilon}(x_1,x_2,x_3)dx_1dx_2dx_3$.

For the type (2) integral, consider typically the domains
$U_+=\{(x_1,x_2,x_3)\in \Bbb E^3|a\le x_1\le b, 0\le x_3\le
\delta, c(x_1)\le x_2\le c(x_1)+d(x_1)x_3\}$ and
$U_-=\{(x_1,x_2,x_3)\in \Bbb E^3|a\le x_1\le b, -\delta\le x_3\le
0, c(x_1)+d(x_1)x_3\le x_2\le c(x_1)\}$ as given in Fig. 7.
\begin{figure}[h]
\begin{center}
\includegraphics[width=0.7\textwidth]{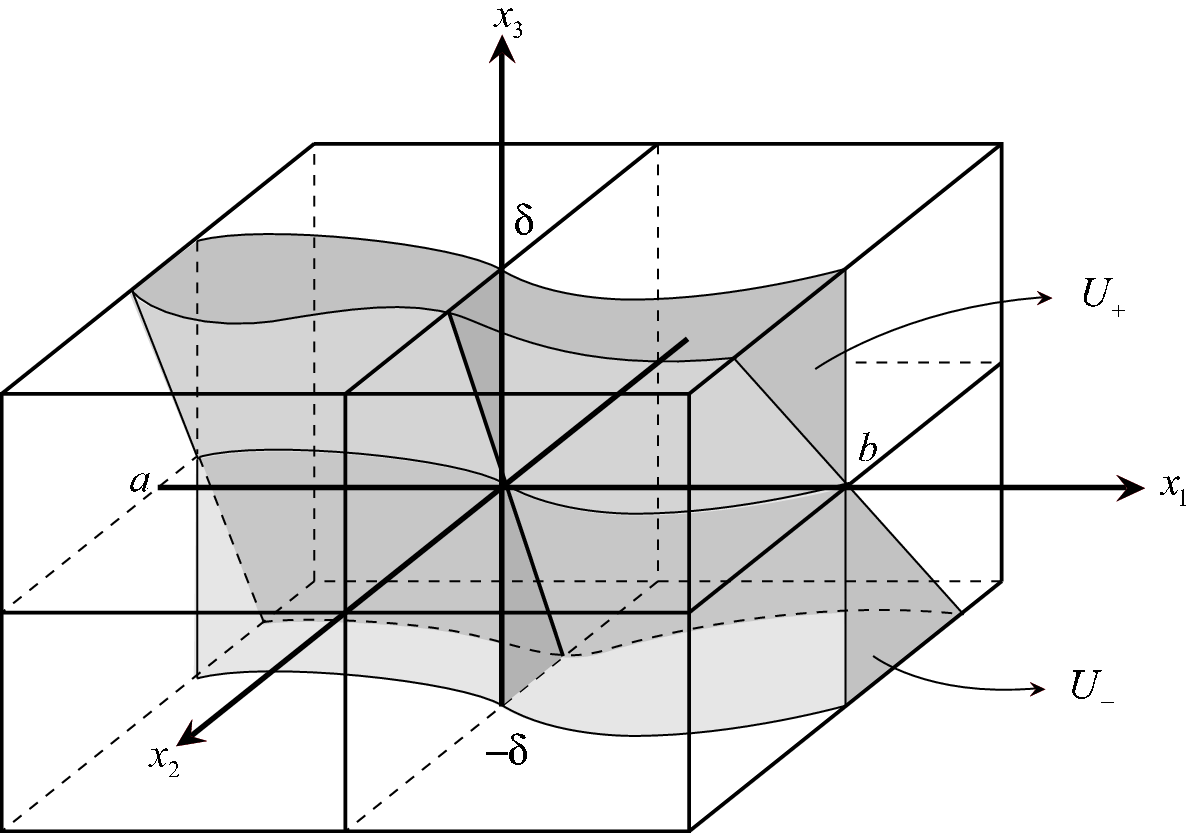}
\caption{\textbf{7}}
\end{center}
\end{figure}

Then the $\text{vol} (U_+)-\text{vol} (U_-)$ is represented by
$$
\int_a^b\!\!\int_0^{\delta}\!\!\int_{c(x_1)}^{c(x_1)+d(x_1)x_3}~d\tilde
V_{\epsilon} -
\int_a^b\!\!\int_{-\delta}^0\!\int^{c(x_1)}_{c(x_1)+d(x_1)x_3}~d\tilde
V_{\epsilon},
$$
and simplified to
\begin{equation}\label{2}
\int_a^b\!\!\int_{-\delta}^{\delta}\!\int_{c(x_1)}^{c(x_1)+d(x_1)x_3}~F_{\epsilon}(x_1,x_2,x_3)~dx_2dx_3dx_1.
\end{equation}
As we have shown in the proof of Step 1 of Theorem 2.1, the pole
of $\int_c^{c+d\cdot x_3}~F_{\epsilon}~dx_2$, as a function of
$x_3$, has negative imaginary part and is analytic on
$\{x_3=\alpha_3+\beta_3 i|\beta_3\ge0\}$,  and hence we have
$$
\int_{-\delta}^{\delta}\!\int_c^{c+d\cdot
x_3}~F_{\epsilon}~dx_2dx_3=\int_{\gamma}\int_c^{c+d\cdot
x_3}~F_{\epsilon}~dx_2dx_3.
$$
From the uniform boundedness of $F_{\epsilon}$, it follows that
$$
\lim_{\epsilon \to 0} \int_a^b\!\!\int_{\gamma}\int_c^{c+d\cdot
x_3}~d\tilde
V_{\epsilon}=\int_a^b\!\!\int_{\gamma}\int_c^{c+d\cdot
x_3}\lim_{\epsilon \to 0}d\tilde
V_{\epsilon}=\int_a^b\!\!\int_{\gamma}\int_c^{c+d\cdot x_3}~d\mu_0
<\infty.
$$

Let us think about the third type integral. In this case we
integrate on the domain lying only in one side $\Bbb H^3$ or
$\Bbb L^3$, and the integral becomes
\begin{equation}\label{3}
\int_a^b\!\!\int_0^{\delta}\!\!\int_{c(x_1)+d(x_1)x_3}^{c(x_1)+d(x_1)
x_3+g(x_1,x_3)}~F_0(x_1,x_2,x_3)~dx_2dx_3dx_1,
\end{equation}
where $g(x_1,0)=\frac{\partial g}{\partial x_3}(x_1,0)=0$ and
$g\in C^{1,\alpha}$ from the hypothesis of
$C^{1,\alpha}$-transversality of the boundary of $U$ and implicit
function theorem for $C^{1,\alpha}$ function. The finiteness of
(3) follows from the finiteness of $\int_0^{\delta}\int_{c+d\cdot
x_3}^{c+d\cdot x_3+g(x_3)}~F_0~dx_2dx_3,$ where
$"g(x_3)"=g(x_1,x_3)$ for each fixed $x_1$ abusing the notation
for $g$. And this integral is equivalent to
$$
\int_0^{\delta}\!\!\int_{c+d\cdot x_3}^{c+d\cdot
x_3+g(x_3)}~\frac{1}{x_3^2}dx_2dx_3=\int_0^{\delta}~\frac{g(x_3)}{x_3^2}dx_3.
$$
The $C^{1,\alpha}$ condition gives us $|g(x_3)|\le C
|x_3|^{1+\alpha}$ and hence
$\int_0^{\delta}~\frac{|g(x_3)|}{x_3^2}dx_3\le\int_0^{\delta}~\frac{1}{x_3^{1-\alpha}}dx_3<\infty,$
and therefore $\int_0^{\delta}~\frac{g(x_3)}{x_3^2}dx_3<\infty$.

We have shown that the local volumes are finite. But this is not
enough in dimension 3. For this type of finitely additive measure
$\mu$ is very subtle and we can not say in general that the volume
of the intersection of two domains with finite volumes is also
finite.

Let's arrange boxes carefully as in the following picture around
the boundary of $U$ and $(x_1,x_2)$-coordinate plane. The picture
is the intersections of boxes with $(x_1,x_2)$-coordinate plane
and shows the wedge shaped domains obtained as intersections
($G_i$) of two boxes and discrepancies ($F_i$) not covered by
boxes.

Notice that $F_i\cap U$, the domain not covered by boxes $S_i$, is
contained in the tetrahedron $T$ which is bounded by the sides of
the boxes and the tangent plane of $\partial U$ (or a suitable
plane so that the tetrahedron $T$ contains $F_i\cap U$). The
domain $G_i\cap U$, overlapped by two boxes, is contained in the
prism minus tetrahedron. We have already shown that the volume of
prism is finite as it is a type (1) integral. Hence if we can show
that the volume of tetrahedron $T$ is finite, then we can complete
the proof of the theorem. But $T$ lies completely in $\Bbb H^3$ or
$\Bbb L^3$ and also being a subset of a cone, it suffices to show
that the cone type domain $E=\{(x_1,x_2,x_3)|0\le x_3\le\delta,
x_3\ge k\sqrt{x_1^2+x_2^2}\}$ has finite volume. Because the
measure in $\Bbb H^n$ or $\Bbb L^n$ is essentially positive
measure.
\begin{equation}\label{4}
\aligned \text{vol
}(E)&=\int_0^{\delta}\!\!\int_0^{\frac{x_3}{k}}\!\!\!\int_0^{2\pi}~\frac{rd\theta
drdx_3}{2x_3^2(r^2+(x_3-1)^2)}\\
    &\sim \int_0^{\delta}\!\!\int_0^{\frac{x_3}{k}}~\frac{r}{x_3^2}drdx_3\\
    &=\int_0^{\delta}~\frac{1}{2k^2}dx_3<\infty
\endaligned
\end{equation}
\end{proof}

\begin{figure}[h]
\begin{center}
\includegraphics[width=0.86\textwidth]{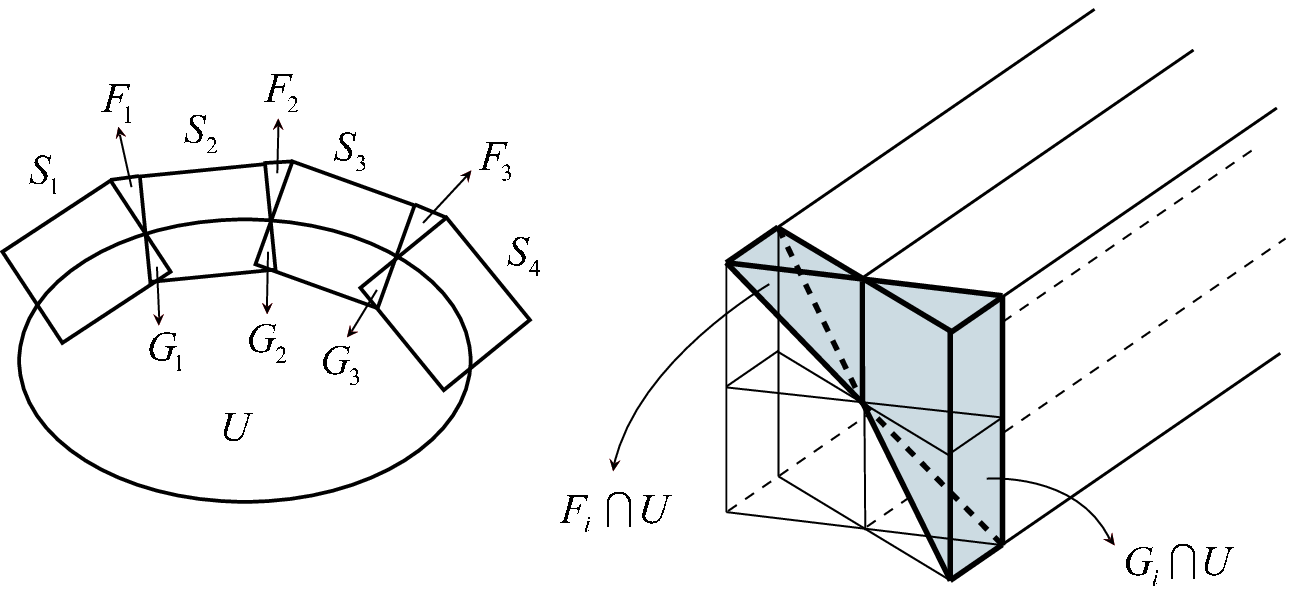}
\caption{\textbf{8}}
\end{center}
\end{figure}

\begin{rem}\label{c} The regularity condition for $\partial U$ is
sharp in the theorem, and in fact there exists a
$C^1$-transversal domain $U$ with infinite volume.

Let $U:=\{(x_1,x_2,x_3)|a\le x_1\le b, -\delta\le x_3\le \delta<1,
-1\le x_2\le h(x_3)\}$, where $h(x_3)=-x_3/\log x_3$ for $x_3\in
(0,\delta]$ and $h(x_3)=0$ for $x_3\in [-\delta,0]$, then it is
easy to show that the volume of $U$ is  infinite by showing that
$$\aligned\mu(U)&\sim\int^b_a\!\!\int^{\delta}_0\!\!\int^{-\frac{x_3}{\log
x_3}}_0~\frac{dx_2 dx_3
dx_1}{x_3^2(x_1^2+x_2^2+(x_3-1)^2)}\sim\int^b_a\!\!\int^{\delta}_0\!\!\int^{-\frac{x_3}{\log
x_3}}_0~\frac{dx_2 dx_3
dx_1}{x_3^2}\\&=\int^b_a\!\!\int^{\delta}_0~\frac{dx_3 dx_1}{-x_3
\log x_3}=\infty.
\endaligned
$$
\end{rem}

\begin{rem}\label{d} For the higher dimensional case, if we can handle a domain $U$ as in the proof Theorem 3.2,
 from the formula of volume form, we may expect the necessary regularity condition for the finiteness as follows:

 $n=4,\quad$ $C^{1,\frac 12+\alpha}$-transversal

$n=5,\quad$ $C^{2,\alpha}$-transversal

$\qquad\quad\qquad\vdots$

The necessary regularity increases by $1/2$ for each dimension
increase. We do not pursue this issue here any more. But note
that this condition is sharp in the sense that we can find a
domain with infinite volume as in Remark \ref{c} if $\alpha=0$.
\end{rem}

\begin{rem}\label{e}
In the proof of Theorem 3.2, we used the volume form $d\tilde
V_{\epsilon}$ in the computation of integrals. But we can use
$d\mu_{\epsilon}$ as well instead of $d\tilde V_{\epsilon}$.
Indeed for the integrals (3) and (4), both $d\tilde V_{\epsilon}$
and $d\mu_{\epsilon}$ will converge to the singular volume form
$d\mu_0$ and gives the same value for the integrals. For the
integrals (1) and (2), the replacement by $d\mu_{\epsilon}$ leads
to the same integral equation by the same pole argument and
uniform boundedness, and eventually get the same integration
value.
\end{rem}

As a final results, we will show that the volume of
$C^{1,\alpha}$-transversal 3-dimensional domain $U$ is invariant
under hyperbolic isometries. Of course we can obtain the same
result in dimension 2 for $C^{0,\frac 12+\alpha}$-transversal
domain $U$ similarly and more easily.

\begin{thm}\label{f} The volume of $C^{1,\alpha}$-transversal domain $U$ is invariant under isometry.
\end{thm}

\begin{proof}
Since the hyperbolic isometries are generated by reflections, we
show the theorem for a reflection $g$. Furthermore it suffices to
show $\text{vol}~(U)=\text{vol}~(gU)$ for each of the four types
domain appeared as (1),(2),(3), and (4) in the proof of Theorem
3.2.

For types (1), we can write as follows:
$$
\aligned \text{vol}~(U)&=\lim_{\epsilon \to
0}\int_B\int^{\delta}_{-\delta}~d\tilde
V_{\epsilon}=\int_B\int_{\gamma}~d\mu_0=\int_B\int_{\gamma}~g^*(d\mu_0)\\
&=\lim_{\epsilon \to 0}\int_B\int^{\delta}_{-\delta}~g^*(d\tilde
V_{\epsilon})=\text{vol}~(gU).
\endaligned
$$
 Here it is enough to give the proof of the
fourth equality, which requires the pole argument and uniform
boundedness as we have used several times before. Indeed notice
that $g^*d\tilde V_{\epsilon}$ never has a real pole for all
reflections $g$ since $g$ is real and $d\tilde V_{\epsilon}$ does
not have a real pole. Hence the poles of $g^*d\tilde
V_{\epsilon}$ have either positive imaginary parts or negative
imaginary parts for all $g$ by continuity with respect to $g$.
Now it suffices to check the sign of imaginary part for a
particular reflection $g_0$ that fixes $(x_2,x_3)$-coordinate
plane. This fixes $x_3$-axis and hence its complexification
preserves negative imaginary part of complex $x_3$-axis  and
poles of $g_0^*d\tilde V_{\epsilon}$ has negative imaginary part
since $d\tilde V_{\epsilon}$ does. The uniform boundedness on a
compact set follows by the same argument as in the proof of
Theorem 2.1.

The invariance of type (2) integral follows similarly.

The domain appeared in the integrals of type (3) and (4) are
either in hyperbolic or Lorentzian space and the integrals are
usual volumes which of course are isometry invariant.

\end{proof}

\begin{flushleft}
{\sc  Yunhi Cho\\Department of Mathematics\\ University of Seoul\\
Seoul 130-743, Korea\\[5pt]}
{\it E-mail}: yhcho@uos.ac.kr\\
\vskip 1pc

{\sc Hyuk Kim\\ Department of Mathematics\\ Seoul National
University\\ Seoul 151-742,
Korea\\[5pt]}
{\it E-mail}: hyukkim@math.snu.ac.kr
\end{flushleft}

\begin{thebibliography}{99}

\bibitem{CK} Yunhi Cho and Hyuk Kim, {\em The Analytic Continuation of Hyperbolic Space},
(preprint).



\end{thebibliography}
\end{document}